# A Parameter Estimation of Fractional Order Grey Model Based on Adaptive Dynamic Cat Swarm Algorithm

Binyan Lin, Fei Gao*, Meng Wang, Yuyao Xiong, Ansheng Li

Department of Mathematics, School of Science, Wuhan University of Technology, Wuhan 430070, China
E-mail: hgaofei@gmail.com

**Abstract:** In this paper, we utilize ADCSO (Adaptive Dynamic Cat Swarm Optimization) to estimate the parameters of Fractional Order Grey Model. The parameters of Fractional Order Grey Model affect the prediction accuracy of the model. In order to solve the problem that general swarm intelligence algorithms easily fall into the local optimum and optimize the accuracy of the model, ADCSO is utilized to reduce the error of the model. Experimental results for the data of container throughput of Wuhan Port and marine capture productions of Zhejiang Province show that the different parameter values affect the prediction results. The parameters estimated by ADCSO make the prediction error of the model smaller and the convergence speed higher, and it is not easy to fall into the local convergence compared with PSO (Particle Swarm Optimization) and LSM (Least Square Method). The feasibility and advantage of ADCSO for the parameter estimation of Fractional Order Grey Model are verified.

## 1 Introduction

Recently, there has been growing interest in the area of parameter estimation of grey system to improve the predictive accuracy[1-4]. Wu, Wang, C.H and Meng, et al respectively have defined fractional order accumulating operators according to Taylor series expansion[5-8], binomial coefficients[9], Gamma function[10], etc. to reduce the perturbation of the forecasting model and promote the stability, and obtain higher accuracy of prediction. The model has been successfully applied in many disciplines to deal with the uncertainty of the system model and incomplete information cases[11-13]. Although considerable research has been devoted to integer order grey model, little attention has been paid to optimize the parameters of the fractional grey system.

There have been many studies utilizing LSM (Least Square Method) to estimate the parameters of the grey model. But LSM aims to minimize the sum of square deviation instead of average relative error, thus the average relative error cannot reach the minimum. Therefore, more scholars have proposed other algorithms to estimate parameters of the integer order grey model to improve the accuracy of model. Zheng, et al. used Genetic Algorithm to estimate the parameters in the case that optimization function is to minimize the average relative error[14]. He, et al. tried linear programming model to estimate parameters[15]. Zhang optimized the background value and boundary value of the grey model and greatly improved the precision of parameters by Particle Swarm Optimization[16]. However, most intelligent algorithms are easily trapped in local optimum.

Cat Swarm Optimization (CSO) is a new swarm intelligence algorithm proposed by Chu，et al[17, 18]. This algorithm has been applied in many fields such as parameter estimation of nonlinear model[19, 20], path planning, image clustering[21, 22] and linear antenna array[23] etc. due to its better local and global searching ability and faster convergence speed.

In this paper, we establish nonlinear objective function and utilize Adaptive Dynamic Cat Swarm Optimization[24] to estimate the parameters of fractional order grey model, which solves the problem that general swarm intelligence algorithms easily fall into the local optimum and optimize the accuracy of the model.

This paper is organized as follows. In Section 2, the fractional order grey model is discussed. In Section 3, the Adaptive Dynamic Cat Swarm Optimization is described. In Section 4, the advantage of utilizing ADCSO to estimate the parameters is clarified by two cases. Finally, Section 5 describes the concluding remarks of this work.

## 2 Grey System Theory and Grey Model

The forecasting model GM (1,1) is one of the most common models in grey system theory. Essentially, GM (1,1) accumulates the original data sequences in order to minimize the influence from random interfering factors, then simulate the trend with exponential curve and find the parameters of the model by least squares estimation or other methods. Therefore, in some way the model can be denoted as an exponential simulating model based on accumulation and least squares method. Compared with traditional GM (1,1), fractional order accumulating grey model is able to reduce perturbation of the prediction model and enhance the stability[5-8], which can provide higher accuracy of prediction.

### 2.1 Discrete grey model based on fractional order accumulation

Let the original sequence be:
$$X^{(0)} = (\xi^{(0)}(1), \xi^{(0)}(2), \cdots, \xi^{(0)}(n)) \quad (1)$$

Then the *r*-order accumulation generating sequence is as follows:

*This work is supported by National Natural Science Foundation (NNSF) of China under Grant 91324201 and the Natural Science Foundation No.2014CFB865 of Hubei Province of China.

Corresponding author Prof. Dr. Fei GAO

$$\xi^{(r)}(k) = \sum_{i=1}^{k} \frac{\Gamma(r+k-i)}{\Gamma(k-i+1)\Gamma(r)} \xi^{(0)}(i) \qquad (2)$$

And the *r*-order reduction generating sequence is as follows:

$$\xi^{(-r)}(k) = \sum_{i=1}^{k-1} (-1)^i \frac{\Gamma(r+1)}{\Gamma(i+1)\Gamma(r-i+1)} \xi^{(0)}(k-i) \qquad (3)$$

Similarly, its nearby mean accumulation generating sequence is:

$$\eta^{(r)}(k) = \frac{1}{2}(\xi^{(r)}(k) + \xi^{(r)}(k-1)) \qquad (4)$$

The *r*-order accumulating grey model has following definition

$$\xi^{(r)}(k) - \xi^{(r)}(k-1) + a\eta^{(r)}(k) = b \qquad (5)$$

as to *r*-order accumulating grey model.

And the whitenization differential equation is:

$$\frac{d\xi^{(r)}(t)}{dt} + a\xi^{(r)}(t) = b \qquad (6)$$

then the time response sequence of *r*-order GM (1,1) is:

$$\hat{\xi}^{(r)}(k+1) = (\xi^{(0)}(1) - \frac{b}{a})e^{-ak} + \frac{b}{a} \qquad (7)$$

where $k = 1, 2, \cdots, n-1$.

The inverse accumulation generating operator is:

$$\hat{\xi}^{(0)}(1) = \xi^{(0)}(1)$$
$$\hat{\xi}^{(0)}(k) = (\hat{\xi}^{(r)})^{-r}(k) = \sum_{i=1}^{k-1}(-1)^i \frac{\Gamma(r+1)}{\Gamma(i+1)\Gamma(r-i+1)} \hat{\xi}^{(r)}(k-i) \qquad (8)$$

where $k = 2, 3, \cdots, n$.[25]

As the parameters $A = [a,b]^T$ of the fractional order grey model directly influence the utility of the model, the estimation of parameter act as an important role in improving the simulating performance of the model. There have been various studies concerning parameter estimation by LSM assuming the data are characterized by normal distribution. But grey model mainly aims to predict small sample data whose characteristics are scarcely distinct and there may be isolated point. As a result, it is not easy to distinguish abnormal value when using least squares estimate, which may finally lead to moderate interference in parameter estimation. For all the reasons above, we estimate the parameters $A = [a,b]^T$ by Adaptive Dynamic Cat Swarm Optimization (ADCSO).

## 3 Adaptive Dynamic Cat Swarm Optimization

Cat Swarm Optimization is a swarm intelligence algorithm based on cat behavior and swarm intelligence[17, 18]. The algorithm has two modes: Seeking mode and tracing mode. And cats are grouped randomly with a certain mixture ratio to find the best position.

In order utilize ADCSO to estimate the parameters of Fractional Order Grey Model, we set up a fitness function describing in (9).

$$f(a,b) = \frac{1}{n-1} \sum_{k=2}^{n} \frac{|\hat{X}_0(k) - X_0(k)|}{X_0(k)} \qquad (9)$$

Then the problem becomes a nonlinear optimization problem, and the goal is to minimize the fitness function.

$$\min f(a,b) = \frac{1}{n-1} \sum_{k=2}^{n} \frac{|\hat{X}_0(k) - X_0(k)|}{X_0(k)} \qquad (10)$$

The two modes of ADCSO are described below:

### 3.1 Seeking mode

Cats have a strong curiosity of moving things, but most of the time they are inactive. They are always alert, and always on guard even when resting[17]. This sub mode is to describe the situation of resting cats which are looking around for next move.

Seeking mode has four general definitions described below:

Seeking memory pool (*SMP*): define the size of seeking memory of each cat, which indicates the solution space of the parameters of Fractional Order Grey Model. The value of *SMP* will not be affected by other parameters.

Seeking range of the selected dimension (*SRD*): define the mutative ratio for the selected dimensions.

Counts of dimensions to be changed (*CDC*): define the counts of dimensions to change, and its value is between 0 to *D*, where *D* is the total dimension.

Self-position considering (*SPC*): *SPC* is a Boolean variable, which decides whether the cat will retain the present position as one of the candidates. [18]

The process of the seeking mode can be described as follow[18, 24]:

Step1: Make *K* copies of the present position of the $cat_i$, $i = 1,2,\cdots,N$, where *K=M* and *M* is the size of seeking memory of $cat_i$. If $SPC_i$ is true, let *K=M-1*, and the present position will be one of the candidates.

Step2: For all positions of seeking memory pool $\{x_i^j\}_{j=1}^M$, $i = 1,2,\cdots,N$, according to the *SRD* and *CDC*, randomly plus or minus $\eta$ (mutative ratio defined by *SRD*) percent of the present values and then obtain the new position by (11),

$$x_i^j \leftarrow x_i^j \pm \eta |x_i^j|, \text{ where } j = 1,2,\cdots,M \qquad (11)$$

Step3: Calculate the fitness of all candidate positions.

Step4: Let each candidate position the same probability, if the fitness of all the candidate positions are exactly equal. Otherwise, calculate the probability of each candidate position by (12)

$$P_j = \frac{|f_{ij} - f_{max}|}{|f_{max} - f_{min}|}, \text{ where } j = 1,2,\cdots,M \qquad (12)$$

Step5: The candidate position with the highest probability is the position the $cat_i$ going to move. Then update the $cat_i$'s old position with new position, and get new value of parameters.

### 3.2 Tracing mode

Although cats are mostly resting, their tracing skills are inherent. Based on that, the other sub mode is to describe the situation of tracing targets and food. In this mode, cats move to the new position quickly[17].

The process of the tracing mode[18] can be described as follow:

Step1: Update the velocities for every dimension by (13)[24].

$$v_{id} \leftarrow \omega_d \times v_{id} + r \times c_d \times (x_{best,d} - x_{id}) \qquad (13)$$

where $d = 1, 2, \cdots, D$, $i = 1, 2, \cdots, N$
$r \sim U[0,1]$.

$v_{id}$ is the $d$th dimension velocity of $cat_i$.

$x_{id}$ is the $d$th dimension position of $cat_i$.

$c_d$ is the adaptive acceleration coefficient and $c_0$ is the initial acceleration coefficient which is usually equal to 2.05.

$$c_d = c_0 - \frac{D-d}{2D}, \text{ where } d = 1, 2, \cdots, D$$

$\omega_d$ is the adaptive inertia weight and $\omega_0$ is the initial inertia weight which is usually equal to 0.6.

$$\omega_d = \omega_0 + \frac{D-d}{2D}, \text{ where } d = 1, 2, \cdots, D$$

Step2: Check if the updated velocity $v_{id}$ exceeds the preset velocity range $[v_{min}, v_{max}]$.

$$v_{id} = \begin{cases} v_{min} & v_{id} \in [-\infty, v_{min}) \\ v_{id} & v_{id} \in [v_{min}, v_{max}] \\ v_{max} & v_{id} \in (v_{max}, +\infty] \end{cases} \quad (14)$$

Step3: Update the position of $cat_i$ by (15)

$$x_{id} \leftarrow x_{id} + v_{id}, \text{ where } d = 1, 2, \cdots, D \quad (15)$$

### 3.3 Core description of ADCSO

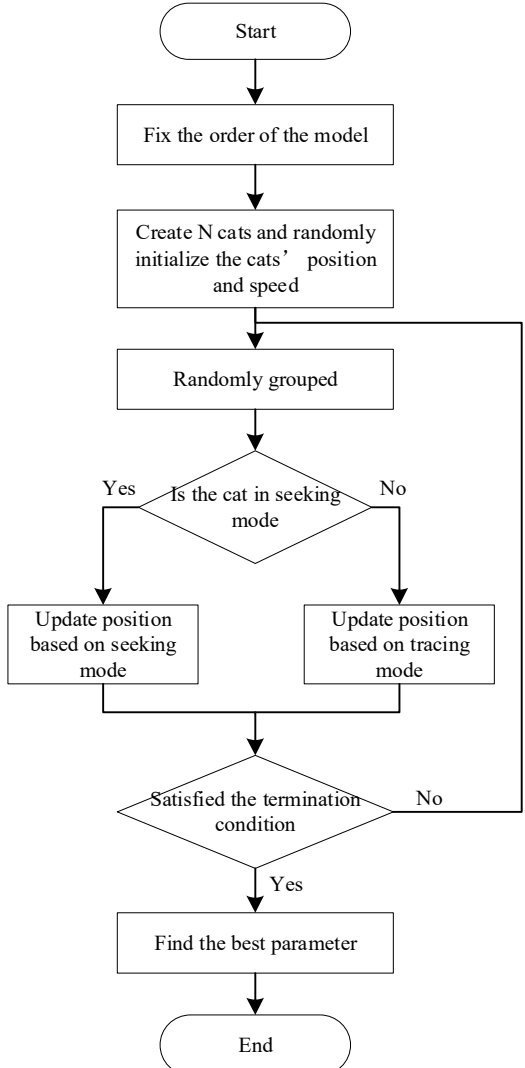

Fig. 1: The process of ADCSO

In order to find the optimal global position, the cats will be divided into seeking group and tracing group according to a certain mixture rate. The cats will update their position according to their group mode and find the optimal solution[18]. The basic process of ADCSO is as follows:

Step1: Create $N$ cats and randomly initialize the cats' positions, velocities and the flags for cats[24]. Then we obtain the initial solution space of the parameters of Fractional Order Grey Model.

Step2: Divide the $N$ cats into two groups (seeking group and tracing group) with a certain mixture ratio (*mr*) randomly.

Step3: Calculate the fitness value of each cat and keep the best cat into memory. And remember the position of the best cat, for it represents the best parameters of Fractional Order Grey Model so far.

Step4: Set the cats into seeking process or tracing process according their flags. If the cat is in tracing mode, set the cat into tracing process, otherwise set it into seeking process.

Step5: Check whether the parameters satisfy the termination condition, if satisfied, terminate the program, otherwise return to Step2. The process is shown on Fig. 1.

## 4 Experiments on numerical optimization

### 4.1 Prediction for Wuhan Port's container throughput

Based on the data about Wuhan Port's container throughput[26] shown in Table 1, we built a fractional order grey model and programmed in Matlab_R2015b, using ADCSO, PSO (particle swarm optimization), and LSM (least square method) respectively for estimating parameters *a* and *b* to calculate the average error of three different orders (0.25, 0.5 and 0.75).

Table 1: Wuhan Port's container throughput

| Year | Container throughput (TEU) |
| --- | --- |
| 2011 | 714700 |
| 2012 | 765000 |
| 2013 | 860412 |
| 2014 | 1005200 |
| 2015 | 1061400 |

Parameter settings for CSO and PSO are shown in Table 2 and Table 3.

Table 2: Parameter settings for CSO

| Parameter | Range |
| --- | --- |
| $N$ | 40 |
| $M$ | 30 |
| $\eta(SRD)$ | 0.2 |
| $mr$ | 0.2 |
| $c$ | 1.05 |
| $w$ | 0.6 |
| $Iter_{max}$ | 300 |

Table 3: Parameter settings for PSO

| Parameter | Range |
|---|---|
| N | 40 |
| c1 | 1.5 |
| c2 | 1.5 |
| w | 0.7 |
| $Iter_{max}$ | 300 |

CSO and PSO are both kinds of swarm intelligence algorithm, so there exist differences among the results of each run when using them. In order to evaluate the estimation of the two algorithms more objectively, we tested each for ten times and take the average of their results. Then put them into comparison, which is shown below.

Table 4: Comparison among the three algorithm's results of data in Table 1

| Algorithm | Error when r=0.25 (%) | Error when r=0.5 (%) | Error when r=0.75 (%) |
|---|---|---|---|
| LSM | 1.57 | 1.36 | 1.47 |
| PSO | 2.03 | 2.52 | 3.10 |
| ADCSO | 1.15 | 1.21 | 1.36 |

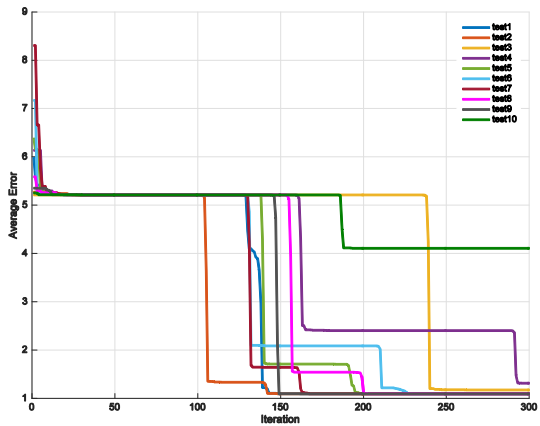

Fig. 2: The iteration figure of PSO of r=0.21

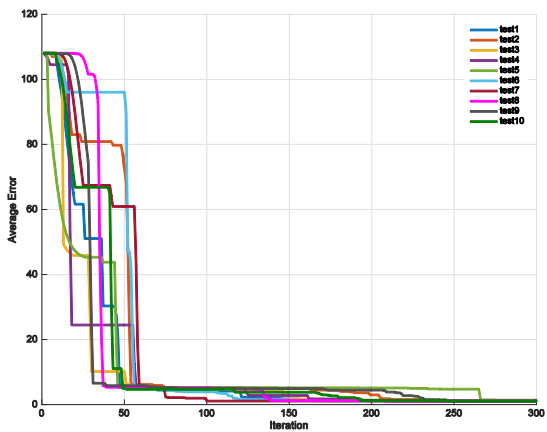

Fig. 3: The iteration figure of CSO of r=0.21.

We can see from Table 4, Fig.2 and Fig.3 that the prediction effect of ADCSO is relatively smaller than others of each order. LSM uses linear approximation to process the data, which leads to a considerable error. PSO sometimes can get decent prediction, but the algorithm is unstable for it often plunges into a local optimum in the calculation process, thus it causes larger average error. Therefore, with higher fitness and faster convergence, CSO is a better algorithm to deal with this type of problem.

Besides, since the order of fractional grey model affects the model performance, we set the step size to 0.01, using ADCSO to optimize the order between 0 and 1, and get the best order $r$=0.21. And the estimated value of parameters for this set of data are

$$r=0.21, a=0.015, b=212927.$$

### 4.2 Prediction for marine capture productions in Zhejiang Province

We also collected the data about marine capture productions in Zhejiang Province during 2007-2013[13], which is shown in Table 5.

Table 5: Marine capture productions in Zhejiang Province during 2007-2013

| Year | Marine capture productions(t) |
|---|---|
| 2007 | 3210300 |
| 2008 | 3272300 |
| 2009 | 3152300 |
| 2010 | 3279100 |
| 2011 | 3411200 |
| 2012 | 3474600 |
| 2013 | 3606700 |

After building the fractional order grey model, we used the three algorithms mentioned above again to make a prediction. The parameter settings are the same as Table 2 and Table 3. And the results are shown below.

Table 6: Comparison among the three algorithm's results of data in Table 5

| Algorithm | Error when r=0.25 (%) | Error when r=0.5 (%) | Error when r=0.75 (%) |
|---|---|---|---|
| LSM | 2.328 | 3.0196 | 2.826 |
| PSO | 2.048 | 3.081 | 4.785 |
| ADCSO | 1.485 | 2.151 | 2.158 |

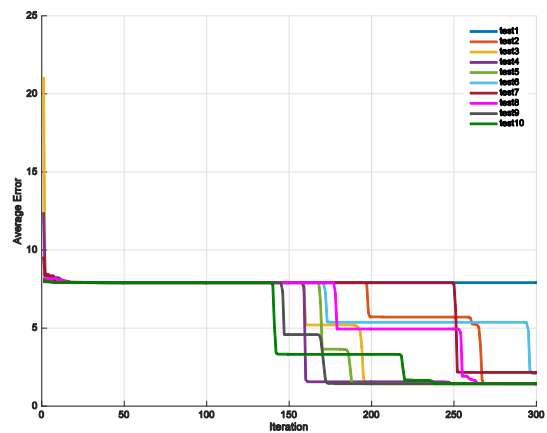

Fig. 4: The iteration figure of PSO of r=0.25

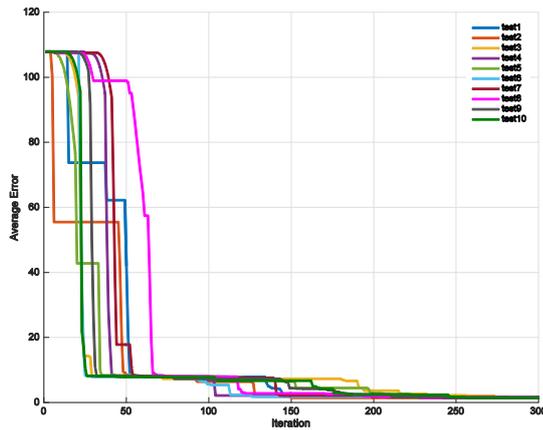

Fig. 5: The iteration figure of CSO of *r*=0.25

By analyzing the information in Table Fig.4 and Fig.5, we can get similar conclusions to the optimization experiment in 4.1. And the estimated value of parameters in this experiment are

$$r=0.06, a=-0.0409, b=3788.$$

From the perspective of the algorithm, particles in PSO updates their own positions according to two extremes—the best position traveled by the individual and of the global traversal, while the tracing mode of CSO is achieved by the global optimal solution, which drive the entire optimization process, speeding up the search. CSO can effectively overcome the shortcomings of other intelligent optimization algorithms (such as inefficient computing, easily into the local optimum, etc.). This is due to the fact that the algorithm has the following advantages: 1. Each solution to the problem is represented by a position of a cat, which is easy to program. A small quantity of parameters has little effect on the algorithm's solving result and efficiency. 2. By dividing all the cats into two groups (seeking mode and tracing mode), ADCSO has better global search capabilities and is less likely to fall into local optimum. And its higher convergence speed determines higher computing efficiency.

## 5 Conclusion

In order to estimate the parameters of fractional order grey model to improve the accuracy of the model, the paper has utilized ADCSO to optimize the parameters. Compared with some algorithms (LSM and PSO), we have found that the parameters estimated by ADCSO are more stable, convergent quickly and can jump out of the local optimum to arrive near optimum solution quickly. The fractional order grey model with optimized parameters can be widely used in small sample prediction in various industries for its higher accuracy. Similarly, ADCSO can be widely used to estimate the parameters of nonlinear systems to improve the accuracy of the systems.